\documentclass[12pt]{article}
\usepackage{amsfonts}
\usepackage{mathrsfs}
\usepackage{amsmath,amssymb}

\baselineskip 5.5pt \lineskip 5.5pt \numberwithin{equation}{section}

\def\QED{\hfill$\Box$\par}
\def\cl{\centerline}
\def\vs{\vspace*}
\def\ni{\noindent}
\def\C{\mathbb{C}}
\def\Z{\mathbb{Z}}
\def\la{\lambda}

\def\bo{{\bf 1}}
\def\rar{\rightarrow}

\def\uu{\mathcal{U}}
\def\pp{\mathcal{P}}

\def\p-{\mathcal{P}_{\leq 0}}
\def\G{\mathfrak{g}}
\def\L{\Lambda}
\def\n{\mathfrak{n}}
\def\h{\mathfrak{h}}
\def\b{\mathfrak{b}}
\def\B{\mathcal{B}}
\def\a{\alpha}

\begin{document}
\cl{{\bf\large WHITTAKER MODULES FOR A LIE ALGEBRA OF BLOCK
TYPE}}\vs{10pt}

\cl{Bin Wang, \,\, Xinyun Zhu}

\begin{abstract}{\footnotesize In this paper, we study Whittaker modules
for a Lie algebras of Block type. We define Whittaker modules and
under some conditions, obtain a one to one correspondence between
the set of isomorphic classes of Whittaker modules over this algebra
and the set of ideals of a polynomial ring, parallel to a result
from the classical setting and the case of the Virasoro algebra.}
\end{abstract}


\noindent{\it{MR(2000) Subject Classification}: 17B10, 17B65,
17B68.}\vs{15pt}

\cl{\bf\S1. Introduction}
\setcounter{section}{1}\setcounter{theo}{0}\setcounter{equation}{0}

Let $\G$ be a Lie algebra that admits a decomposition $$\G = \b_-
\oplus \n $$where $\b_-, \n$ are two Lie subalgebras. Let $$\varphi
: \n \rar \C$$be a homomorphism of Lie algebras. For a $\G$ module
$V$ and $v\in V$, one says that $v$ is a Whittaker vector of type
$\varphi$ if $\n$ acts on $v$ through $\varphi$. A Whittaker module
is then defined to be a module generated by a Whittaker vector.

The category of Whittaker modules for a given algebra (say, $\G$)
admits an initial object and we call it a universal Whittaker module
(say, $M$). Then $M$ is isomorphic to $\uu (\b_-)$ as $\b_-$
modules. Let $Z$ stand for the center of $\G$, and $Z'' = Z \cap
\b_-$. Suppose $\G$ possesses the following properties:

1) for each ideal $I$ of $S(Z'')$, then every Whittaker vector of
the Whittaker module $M/IM$ is of form $p\bar{w}$, with $p \in
S(Z'')$, where $w$ is a Whittaker generator of M;

2) for each $I \subset S(Z'')$, then any nontrivial submodule of
$M/IM$ admits a nonzero Whittaker vector.\vspace*{5pt}

\ni Then it is not hard to set up a correspondence between the set
of isomorphic classes of Whittaker modules and the one of all the
ideals of $S(Z'')$.

In this paper, we consider a Lie algebras of Block type,
$\mathcal{B}$, which is an infinite-dimensional Lie algebra with a
basis $\{x_{a,i}\,|\,a\in\Z,i\in\mathbb{N} \}$ and brackets
\begin{eqnarray}
&&[x_{(a,i)},x_{(b,j)}]=\big((b - 1)i - (a - 1)j\big)x_{(a+b,
i+j-1)}.\label{b1}
\end{eqnarray}
and it has the following decomposition
\begin{eqnarray}
\mathcal{B}=\n_-\oplus \h \oplus \n,
\end{eqnarray}
where
\begin{eqnarray*}
&&\h=\mbox{Span}_\C\{x_{(a,i)},\,| \, a+i = 1 \},\\
&&\n=\mbox{Span}_\C\{x_{(a, i)} \,|\, a+i > 1 \},\\
&&\n_-=\mbox{Span}_\C\{x_{(a, i)} \,|\, a+i < 1 \}.
\end{eqnarray*}

Suppose $\varphi$ is a given good character ( for its definition,
see 2.3.2). Our main result is to show that $\mathcal{B}$ satisfies
the above properties 1) and 2) for such $\varphi$ (see \S3), and
hence obtain a correspondence between Whittaker modules and ideals
of a polynomial ring (of $x_{1,0}$). This is treated in \S4. Note
that for general characters, property 1) may not hold.

Whittaker modules were first discovered for $\mathfrak{sl}_2{(\C )}$
by Arnal and Pinzcon in \cite{AP}. Block showed, in \cite{B} that
the simple modules for $\mathfrak{sl}_2(\C )$ consist of highest
(lowest) weight modules, Whittaker modules and a third family
obtained by localization. This illustrates the prominent role played
by Whittaker modules.

Kostant defined Whittaker modules for an arbitrary
finite-dimensional complex semi-simple Lie algebra $\mathfrak{g}$ in
\cite{K}, and showed that these modules, up to isomorphism, are in
bijective correspondence with ideals of the center
$Z(\mathfrak{g})$. In particular, irreducible Whittaker modules
correspond to maximal ideals of $Z(\mathfrak{g})$. In the quantum
setting, Whittaker modules have been studied by Sevoystanov for
$\uu_h(\mathfrak{g})$ \cite{S} and by M. Ondrus for
$U_q(\mathfrak{sl}_2)$ in \cite{O}. Recently Whittaker modules have
also been studied by M. Ondrus and E. Wiesner for the Virasoro
algebra in \cite{OW}, X. Zhang and S. Tan for
Schr\"{o}dinger-Virasoro algebra in \cite{ZT}, K. Christodoulopoulou
for Heisenberg algebras in \cite{C}, and by G. Benkart and M. Ondrus
for generalized Weyl algebras in \cite{BO}.

We note that our proofs differ from the ones in the classical
setting in the use of the center of the universal enveloping
algebra. The reasoning for this is similar to the one explained in
\cite{OW}. Also, our approach to obtaining property 2) is same as in
\cite{Wang}, different from \cite{OW}.

The paper is organized in the following way. In section 2, we define
Whittaker vectors and Whittaker modules for a class of Lie algebras,
and also construct a universal Whittaker module for them. Then the
Whittaker vectors in a universal Whittaker module are examined in
section 3 and the irreducible Whittaker modules are classified in
section 4. In the last section we discuss some examples.


 \vs{18pt}

\cl{\bf\S2. Preliminaries}

\ni 2.1. {\bf Q-graded Lie algebras} \vspace{6pt}

\ni 2.1.1. Let $V$ be a vector space over $\C$ and $Q$ a free
abelian additive semigroup. By a $Q$-grading of $V$ we will
understand a family $\{V_\a | \a\in Q \}$ of subspaces of $V$ such
that $V = \oplus_{\a\in Q}V_{\a}$. For a nonzero vector $v \in
V_\a$, we say $v$ is a homogeneous vector of degree $\a$. Let $\G$
be a Lie algebra over $\C$ and let $\{ \G_\a \,|\, \a\in Q \}$ be a
grading of $\G$ (as a vector space). Call $\G$ a $Q$-graded Lie
algebra if $[ \G_\a, \G_\beta ] \subset \G_{\a +\beta}$, for all
$\a, \beta \in Q$.

Now suppose $Q$ is totally ordered abelian group by the ordering
$\leq$ that is compatible with its additive group structure. Given a
$Q$-graded Lie algebra, $\G = \oplus_{\a\in Q}\G_{\a}$, and  a
homomorphism of abelian groups $\pi : Q \rar \Z$ that preserves the
ordering, write $\G_m = \sum\limits_{\pi (\a )=m} \G_{\a} $. Then
the Lie algebra $\G = \oplus_{i\in \Z}\G_i$ can be viewed as a
$\Z$-graded Lie algebra too.

\ni 2.1.2.  We now consider the algebra $\mathcal{B}$ defined by
$\eqref{b1}$. Let $Q = \Z \times \Z$, and $\pi : Q \rar \Z$ by $\pi
((a,i))= a + i -1$. Equip $Q$ with a group structure by
$$ (a, i) * (b, j) = (a+b, i+j-1)$$and a ordering by
$$ (a,i) < (b,j) \,\, \mbox{if either}\,\, a+i < b+j
\,\,\mbox{or} \,\, a+i=b+j, i<j .$$ Then $Q$ becomes a totally
ordered abelian group, and $\pi$ preserves the ordering and group
structure. Let $Q' = \{\a \in Q \,|\, \pi(\a )>0 \}, Q'' = \{\a \in
Q \,|\, \pi (\a)\leq 0 \}$ and $K_n = \{ (a,i) \in Q \,|\, i\geq n
\}, n\geq 0$. Set $K'_n = K_n \cap Q', K''_n = K_n \cap Q''$ and $R
= \{(1,0)\}$. So $\B$ (resp. $\n$, resp. $\b_-$)  is a $Q$-graded
(resp. $Q'$-graded, resp. $Q''$-graded) Lie algebra. Write $K=K_0,
K'=K'_0, K''=K''_0.$

\ni 2.2. {\bf{Partitions}}.\vspace{6pt}

\ni 2.2.1. Let $\L$ be a totally ordered set. We define a partition
of $\L$ to be a non-decreasing sequence of elements of $\L$,
$$\mu=(\mu_1,\mu_2,\cdots,\mu_r),\,\, \mu_1\leqslant\mu_2\leqslant\cdots\leqslant\mu_r.$$
Denote by $\pp (\L )$ the set of all partitions. For $\lambda
=(\la_1,\cdots ,\la_r) \in \pp (\L )$, we define the length of $\la$
to be $r$, denoted by $\ell (\la)$, and  for $\a\in \L$, let $\la
(\a)$ denote the number of times $\a$ appears in the partition.
Clearly any partition $\la$ is completely determined by the values
$\la (\a), \a\in \L$. If all $\la (\a)=0$, call $\la$ the null
partition, denoted by $\bar{0}$. Note that $\bar{0}$ is the only
partition of length $=0$. We consider $\bar{0}$ an element of
$\pp(\L )$.

Back to the situation of a Lie algebra $\G$. Define the symbols
$x_{\la}$, for all partitions. For $\bar{0} \neq \la$, define
$x_{\la}$ to be an element of $\uu ( \G )$, the universal enveloping
algebra of $\G$, by
$$x_{\la}=x_{\la_1}x_{\la_2}\cdots x_{\la_r}=\prod\limits_{\a \in
K}x_\a^{\la(\a)} \in \uu(\G)$$ whenever each $x_{\la_i}$ is well
understood as an element of $\G$. And let $x_{\bar{0}}=1 \in \uu{(\G
)}$.

 By PBW theorem, we know that $\{x_{\la} \,|\, \la \in
\pp(K\setminus R) \}$
 form a basis of $\uu (\B )$ over $S(Z)$ where $Z=\C x_{(1,0)}$ and
 $S(Z)$ is the polynomial ring of $x_{(1,0)}$. \vspace*{6pt}

\ni 2.2.2.\, $\uu{(\G)}$ (denoted by $\uu$) naturally inherits a
grading from the one of $\G$. Namely, for any $\a \in Q$, set
$\uu_\a =Span_{\C}\{x_1x_2\cdots x_k \,|\, x_i\in \G_{\a_i}, 1\leq
i\leq k, \sum\limits_{i=1}^{k}\a_i = \a \}$, and then $\uu =
\bigoplus\limits_{ \a \in Q} \uu_\a$ is a $Q$-graded algebra, i.e.
$\uu_\a \uu_\beta \subseteq \uu_{\a + \beta}$. Similarly, $\uu (\n
)$ (resp. $\uu ({\b_-})$) inherit a grading from $\n$, (resp. $\b_-$
). If $x \in \uu_\a$, then we say $x$ is a homogeneous element of
degree $\a$. Set $|\bar{0}|=0$ and $|\la |=\la_1 + \la_2 +\cdots +
\la_{{\ell}(\la)}, \, \forall \la\neq \bar{0}$. Then $x_\la$ is a
homogeneous element of degree $|\la |$. If $u (\neq 0 )$ is not
homogeneous but a sum of finitely many nonzero homogeneous elements,
then denote by $mindeg(u)$ the minimum degree of its nonzero
homogeneous components.

Now let us, for convenience, call any product of elements $x_\a^s$ (
$\a \in K\setminus R, s\geq 0$) in $\uu$ and elements of $S(Z)$ a
monomial, of height equal to the sum of the various $s$'s occurring.
Then we have, by PBW theorem,\vspace{5pt}

\ni {\bf Lemma}. \ \ {\it For $\a, \beta \in K\setminus R, t,k\geq
0$ , $x_\beta^t x_\a^k$ is a $S(Z)$-linear combination of $x_\a^k
x_{\beta}^t$ along with other monomials of height $< t+k$.} \QED

This allows us to make the following definition. If $ x \in \uu $ is
a sum of monomials of height $\leq l$, we say $ht(x)\leq
l$.\vspace{6pt}

\ni 2.2.3. We need some more notation. For $\la =(\la_1,\la_2,\cdots
\la_r)\in \pp (K) , 0< i\leq r, 0\leq j < r, $write
\begin{eqnarray*}
&&\la \{i\}=(\la_1 ,\cdots , \la_i ), \, \la \{0\}=\bar{0} \\
&&\la [j]=(\la_{j+1},\cdots ,\la_r ), \, \la [r]=\bar{0}\\
&&\la <i>=(\la_1,\cdots ,\la_{i-1},\hat{\la_i},\la_{i+1},\cdots
\la_r ).
\end{eqnarray*}\vspace{5pt}

\ni {\bf Lemma }\ \  {\it  Write $\uu''$ for $\uu (\b_-)$. Let
$0\neq x \in \n_\beta, 0\neq y \in \uu''_\gamma$ with $\pi (\beta)
> 0, \pi(\gamma) \leq 0$.

1) if $ s=\pi (\beta + \gamma ) > 0 $, then $[x, y] = \sum\limits_{
s\leq \pi(\a )\leq \pi (\beta) } u_\a $ with $u_\a =
\sum\limits_{i}v^{(\a,i)}w^{(\a,i)}$ where $w^{(\a,i)} \in \n_\a,
v^{(\a,i)} \in \uu'' _{\beta + \gamma - \a}$ and $ht(v^{(\a,i)}) <
ht(y)$ if $ v^{(\a,i)} \neq 0$;

2)if $ \pi (\beta + \gamma ) \leq 0 $, then $[x, y] = \sum\limits_{
0< \pi(\a )\leq \pi (\beta) } u_\a  + u$ with $u \in \uu''_{\beta +
\gamma}$ and $u_\a = \sum\limits_{i}v^{(\a,i)}w^{(\a,i)}$ where
$w^{(\a,i)} \in \n_\a, v^{(\a,i)} \in \uu'' _{\beta + \gamma - \a}$
and $ht(v^{(\a,i)}) < ht(y)$ if $ v^{(\a,i)} \neq 0$}.

\ni {\bf Proof}\ \  Write $y = \sum_{\la}f_{\la}x_\la$ where
$f_{\la} \in S(Z), \la \in \pp(K''\setminus R)$.

But for any $x \in \n_\beta$, one has,
\begin{eqnarray*}
[x, x_\la ] &=& \sum_{i=1}^{\ell (\la )}x_{\la\{i-1\}}[x,
x_{\la_i}]x_{\la[i]}\\
&=& \sum_{i=1}^{\ell (\la )}x_{\la <i-1>}[x, x_{\la_i}] +
\sum_{i=1}^{\ell (\la )}x_{\la\{i-1\}} [[x, x_{\la_i}], x_{\la[i]}]
\end{eqnarray*}
Then one can easily deduce that the lemma follows. \QED
\vspace{10pt}

\ni 2.3 {\bf{Whittaker Module}} \ \

\ni 2.3.1. {\bf{Definition}} \ \ {\it  Given a Lie algebra
homomorphism $\varphi: \n \rar \C$, for a $\B$-module $V$, a vector
$v\in V$ is said to be a Whittaker vector of type $\varphi$ if
$xv=\varphi(x)v$ for all $x\in \n$. Furthermore, if $v$ generates
$V$, then we call $V$ a Whittaker module of type $\varphi$ and $v$ a
cyclic Whittaker vector of $V$.}\vspace{6pt}

\ni 2.3.2. One says a Lie algebra homomorphism $\varphi: \n \rar \C$
is nonsingular, if $\varphi(x_{(a,i)})\neq 0$, for all $(a,i)\in Q$
with $a+i = 2$. For $n\geq 1, s\geq 1$, let $\varphi_{m,j}^{(n,s)}=
\varphi (x_{(2-a,a)})$, where $a=m+j+s-1$, for $m\geq 1, 0\leq j\leq
n$. Denote by $H^{(n,s)}$ the $\infty \times (n+1)$ matrix whose
$(m,j)$ entry is $\varphi_{m,j}^{(n,s)}$. \vspace*{6pt}

{\bf{Definition}} \ \ {\it A Lie algebra homomorphism $\varphi: \n
\rar \C$ is said to be a good character if $\varphi$ is nonsingular
and for all $n\geq 1, s\geq 1$, rank$(H^{(n,s)})=n+1$.}

For a given $\varphi:\n \rar\C$, define $\C_\varphi$ to be the
one-dimensional $\n$-module given by the action $xa= \varphi(x) a$
for all $x\in \n$ and $a\in\C$. Then the induced $\B$-module,
\begin{eqnarray*}
M_\varphi=\uu(\B)\otimes_{\uu(\n)}\C_\varphi,
\end{eqnarray*}
is a Whittaker module of type $\varphi$ with the cyclic Whittaker
vector $w=\bo\otimes 1$.  By PBW theorem, it's easy to see that
$\{x_{\la}w \,|\,\la \in \pp (K_0''\setminus R) \}$ is a basis of
$M_\varphi$ over $S(Z)$ where $Z=\C x_{(1,0)}.$

Besides, for any  ideal $I$ of $S(Z)$, define $L_{\varphi, I} =
M_\varphi/IM_\varphi$ and denote  by $p_I$ the canonical
homomorphism. Then $L_{\varphi, I}$ is a Whittaker module for $\B$.
The following lemma makes $M_{\varphi}$ become a universal Whittaker
module. \vspace{6pt}

\ni {\bf Lemma} {\it Fix $\varphi$ and $M_\varphi$ as above. Let $V$
be a Whittaker module of type $\varphi$ generated by a Whittaker
vector $w'$. Then there is a unique map $\phi:M_\varphi\rightarrow
V$ taking $w=1\otimes 1$ to $w'$.}

\ni{\bf Proof}.\ \ Uniqueness is obvious. Consider $u\in \uu (\B
)$. One can write, by PBW,
$$u=\sum\limits_\alpha b_\alpha n_\alpha,\,\, b_\alpha\in \uu (\b_- ), n_\alpha\in \uu(\n)$$
If $uw=0$, then $uw=\sum\limits_\alpha b_\alpha
\varphi(n_\alpha)w=0$, and therefore $\sum\limits_\alpha b_\alpha
\varphi(n_\alpha)=0$. Now it's easy to see that the map $\phi :
M_{\varphi} \rar V$, defined by $\phi(uw) = uw'$, is well
defined.\QED\vspace{6pt}

\ni 2.3.3. Let $A=S(Z), Z=\C x_{(1,0)}$ and $I$ be an ideal of $A$.
Write $M = M_{\varphi}, \p- = \pp (K''\setminus R)$ and $w' = p_I w
\in M/IM.$\vspace{6pt}

\ni {\bf Lemma} \ \ {\it $M/IM$ admits a basis, $\{ x_\la w'\,|\,
\la \in \p- \},$ over $A/I$.}

\ni{\bf Proof} \ \ Note that $M = \bigoplus\limits_{\la \in
\p-}Ax_\la w$. Hence,
$$M/IM = A/I \otimes_A M = A/I \otimes_A (\bigoplus_{\la \in
\p-}Ax_\la w)=\bigoplus_{\la \in \p-}(A/I) x_\la w' .$$ Then the
lemma follows immediately. \QED\vspace{6pt}

\ni 2.3.4. Assume now that $I$ is an ideal of $A = S(Z)$. Write
$\uu'$ for $\uu (\n)$, and $\uu''$ for $\uu (\b_-)$. Then $V = M/IM$
has a natural grading as a vector space. Namely, based on Lemma
2.3.2, let, for any $\a \in Q''$, $V_\a = \{ x = \sum\limits_{\la
\in \p-}a_\la x_\la w \,|\, a_\la \in (A/I), a_\la = 0
\,\,\mbox{if}\,\,|\la | \neq \a \}$ and then clearly $V
=\bigoplus\limits_{\a\in Q''} V_{\a}$. We say that a nonzero
homogeneous vector $v$ in $M$ is of degree $\a$ if $v \in V_\a$. If
$v\, (\neq 0 )$ is not homogeneous but a sum of finitely many
nonzero homogeneous vectors, then define $mindeg(v)$ to be the
minimum degree of its nonzero homogeneous components. Meanwhile, for
any nonzero vector $v \in V$, let $d(v) = mindeg(v)$ and then there
uniquely exist $v_i \in \uu''_\a,  \a \geq d(v)$ such that $v =
\sum\limits_{d(v)\leq \a} v_{\a}w$ with $v_{d(v)} \neq 0$. Then
define $\ell (v) = ht( v_{d(v)})$. Note $V$ can also be equipped
with a $\Z$-grading through $\pi : Q \rar \Z$. With this grading, we
can introduce notation $mindeg_1 (v)$ and $\ell_1 (v)$, for each
$0\neq v \in V$, parallel to $mindeg(v)$ and $\ell (v)$
respectively.

\vs{18pt}

\cl{\bf\S3. Whittaker Vectors in $M_{\varphi}$ and $L_{\varphi, I}$}
\setcounter{section}{3}\setcounter{theo}{0}\setcounter{equation}{0}

In this section, we characterize the Whittaker vectors in Whittaker
modules for $\B$,  where $\varphi$ is a fixed nonsingular Lie
algebra homomorphism from $\n \rar \C$.  Let $M = M_{\varphi},
w=\bo\otimes 1 $, and $Z=\mathfrak{Z}(\B )=\C x_{(1,0)}$, the center
of $\B$. Notation as in 2.1.2. \vspace{10pt}

\ni 3.1. Assume that $I$ is an ideal of $A = S(Z)$. Set $V = M/IM$
and $w' = p_I ( w )$. Write $\pp_1 = \pp (K_1''), \pp_2 = \pp
(\bar{K}\setminus R)$ where $\bar{K}=K''\setminus K''_{1}$. Then we
have the following lemma. \vspace{6pt}

\ni{\bf Lemma} \ \ {\it   Assume $\varphi$ is a good character, then
every Whittaker vector of $V$ is of form $pw'$ with $p \in A =
S(Z)$.}

\ni{\bf Proof} \ \ Suppose $w''$ is a Whittaker vector of $V$. We
can write, by Lemma 2.3.3, \begin{eqnarray} w'' = \sum_{\la \in
\pp_1,\mu\in \pp_2 }p_{\la, \mu}x_{\la}x_\mu w' \label{3.1}
\end{eqnarray}
where $ p_{\la, \mu} \in A/I$. Obviously it is enough to show that
$p_{\la, \mu} = 0$ if either $\la\neq \bar{0}$ or $\mu \neq
\bar{0}$.

Case a), suppose there exists a $\mu \neq \bar{0} $ such that
$p_{\la,\mu}\neq 0$ for some $\la \in \pp_1.$

Let $\L_1 =\{ (a,i) \,|\, a+i\leq 0, i\geq 1\}, \L_2 =\{ (a,i) \,|\,
a+i\leq 1, i\geq 2\}$ and $\L = \L_1 \cup \L_2 .$ Then $K''_1=\L
\cup \{(0,1)\}$. Set $\pp'_1 = \pp(\L)$. Then obviously one can
rewrite equation (3.1) as
\begin{eqnarray} w'' = \sum_{\la \in
\pp'_1,\mu\in \pp_2 }\sum_{s\geq 0}q_{\la, \mu, s}x_{\la}y^sx_\mu w'
\end{eqnarray}
where $y=x_{(0,1)}, q_{\la, \mu, s}\in A/I$. Then there exists a
$\mu\neq 0$ such that $q_{\la, \mu, s}\neq 0$ for some $\la, s$ by
our assumption.

Now take a $N > 2$ so that for any $\la \in \pp'_1$, if $\exists \,
\la_{i} = (a,j)$ with $j > N-2$, then $q_{\la, \mu, s} = 0 .$
Obviously this can be achieved since there are only finitely many
nonzero $q_{\la, \mu, s}$. Put $u= x_{(2-N, N)} .$ Consider
\begin{eqnarray*}
(u - \varphi(u))w'' =& \sum &q_{\la,\mu,s}[u, x_\la y^sx_\mu ]w'\\
=&\sum & q_{\la,\mu,s}[u, x_\la ]y^sx_\mu w'\\
&+&\sum q_{\la,\mu,s}x_\la [u,y^s ]x_\mu w'\\
&+&\sum q_{\la,\mu,s}x_\la y^s[u, x_\mu ]w' .
\end{eqnarray*}
But note that it can be easily showed by induction that $[u, y^s ] =
f_s(y)u$, where $f_s(y)$ is a polynomial of $y$ with
$deg(f_s(y))\leq s-1.$ Therefore,
\begin{eqnarray}
(u - \varphi(u))w'' = &\sum & q_{\la,\mu,s}[u, x_\la ]y^sx_\mu w'\label{f1}\\
&+&\sum q_{\la,\mu,s}x_\la f_s(y)x_\mu\varphi(u)w' + \sum q_{\la,\mu,s}x_\la f_s(y)[u, x_\mu ] w'\label{f2}\\
&+&\sum q_{\la,\mu,s}x_\la y^s[u, x_\mu ]w' .\label{f3}
\end{eqnarray}
Let $\pi_1 : Q \rar \Z$ by $\pi_1(a, i) = a.$ Set
$$t = min\{\pi_1(\mu_1)\,|\, \mu \in \pp_2, \mbox{and} \,\,
q_{\la,\mu,s}\neq 0, \,\mbox{for some}\,\, \la,s . \}$$ Then $t\leq
0.$ Take a $\tau \in \pp_2$ such that $\tau_1 = (t, 0)$ and
$q_{\la,\mu,s}\neq 0$ for some $\la, s.$ Let $r$ be the maximum
integer such that there exists a $\la$, s.t. $q_{\la,\tau, r} \neq
0$. Set $\tau' = \tau<1> = (\tau_2,\cdot\cdot\cdot
\tau_{\ell(\tau)})$ and $\a = (t-N+2, N-1).$

Note that $\{ x_\la y^kx_\mu \,|\, \la\in \pp'_1, \mu\in \pp_2,
s\geq 0\}$ form a basis of $V = M/IM$ over $A/I$. Clearly, under
this basis, the representation of the formula $\eqref{f3}$ contains
nonzero terms involving $x_\a y^rx_{\tau'}$ which are linearly
independent from each other. However, it is easy to see that for
$\eqref{f1}$ and $\eqref{f2}$, there are no terms involving $x_\a
y^rx_{\tau'}$. Hence, $(u - \varphi(u))w'' \neq 0.$ This contradicts
with the assumption that $w''$ is a Whittaker vector.

Case b), suppose $p_{\la, \mu}= 0$ whenever $\mu \neq \bar{0}$, and
there exists at least a $\la \neq \bar{0}$ such that $p_{\la,
\bar{0}}\neq 0.$

Let $L=\{(a,i) \,|\, a+i=1, i\geq 1\}$ and $L'=\{(a,i) \,|\, a+i <
1, i\geq 1\}$. Set $\mathcal{Q} = \pp (L)$ and $\mathcal{Q}' = \pp
(L')$. Then one can rewrite equation (3.1) as
\begin{eqnarray*}w'' =
\sum_{\la \in \mathcal{Q}',\mu \in \mathcal{Q} }f_{\la,
\mu}x_{\la}x_{\mu} w'
\end{eqnarray*}
 where $f_{\la, \mu}\in A/I$.

i). Assume that there exists a $\la(\neq \bar{0}) \in \mathcal{Q}'$
such that $f_{\la, \mu}\neq 0$ for some $\mu \in \mathcal{Q}$.

Define $\pi_2: Q\rar \Z$ by $\pi_2(a,i)=i$. Take a $n_0 > 0$ such
that $f_{\la, \mu} = 0$ if either there is a $i$ such that $n_0 \leq
\pi_2(\la_i)$ or there is a $j$ such that $n_0 \leq \pi_2(\mu_j)$.
Let $\a_0 = max \{\la_{\ell(\la)} \,|\, \exists \, \mu s.t. f_{\la,
\mu}\neq 0 \}, \tau_0 = (2-n_0, n_0)$, and $y=x_{\tau_0}$. Consider
\begin{eqnarray}(y-\varphi(y))w'' = &&\sum_{\la,\mu} f_{\la,\mu}[y, x_\la x_\mu] w' \nonumber \\
=&&\sum f_{\la,\mu}[y, x_\la]x_\mu w' \label{g1}\\
&&+\sum f_{\la,\mu}x_\la [y, x_\mu ]w' .\label{g2}
\end{eqnarray}

Note that $\triangle = \{ x_\la x_\mu x_\gamma \,|\, \la\in
\mathcal{Q}', \mu\in \mathcal{Q}, \gamma \in pp_2\}$ form a basis of
$V = M/IM$ over $A/I$. Clearly, under this basis, the representation
of the formula $\eqref{g1}$ contains nonzero terms involving
$x_{\tau_0*\a_0}$ which are linearly independent from each other.
However, it is easy to see that for $\eqref{g2}$, there are no terms
involving $x_{\tau_0*\a_0}$. Hence, $(u - \varphi(u))w'' \neq 0.$
This contradicts with the assumption that $w''$ is a Whittaker
vector.

ii). Assume $f_{\la, \mu} = 0$ if $\la \neq \bar{0}$, and there
exists at least a $\mu \neq \bar{0}$ such that $f_{\bar{0}, \mu}\neq
0$.

In this case, we write
\begin{eqnarray*}w'' =
\sum_{\mu \in \mathcal{Q} }b_{\mu}x_{\mu} w'
\end{eqnarray*}
where $b_\mu = f_{\bar{0}, \mu}$.

Let $\sigma_0 = (1-s,s)$ for some $s\geq 1$, and $\sigma_0' =
(1-s-n, s+n)$ for some $n\geq 1$ be such that for all $\mu \in
\mathcal{Q}$, if either $\mu_1 < \sigma_0$ or $\mu_{\ell(\mu)} >
\sigma_0'$, then $b_\mu =0$. Set $y_m = x_{(2-m, m)}, m\geq 1$.
Consider
\begin{eqnarray}(y_m - \varphi(y_m))w'' = &&\sum_{\mu} b_{\mu}[y_m, x_\mu] w' \nonumber \\
=&&\sum_{\mu}\sum_i^{\ell(\mu)}b_{\mu}x_{\mu \{i-1\}}[y_m , x_{\mu_i} ]x_{\mu [i]}w' \nonumber\\
=&&\sum_{\mu}\sum_i^{\ell(\mu)}b_{\mu}x_{\mu <i>}\varphi ([y_m ,
x_{\mu_i}])w' \label{h1}\\
 &&+\sum_{\mu}\sum_{i=1}^{\ell(\mu)}b_{\mu}x_{\mu \{i-1\}}[[y_m , x_{\mu_i} ], x_{\mu
 [i]}]w'. \label{h2}
\end{eqnarray}

 Let $l = max \{\ell(\mu) \,|\, b_\mu \neq 0 \}$ and
 take a $\la \in \mathcal{Q}$ such that $\ell(\la)=l-1$, and
 $\exists \, \mu$ \,s.t.\, $b_\mu \neq 0, \la = \mu <i>$ for some $i$. Now
 set $\sigma_i = (1-s-i, s+i), 0\leq i \leq n$ and let $0\leq i_1 \leq i_2 \leq \cdots \leq i_k$ be such that
 $\{\sigma_{i_j} \,|\, 1\leq j \leq k \} = \{\a \,|\, \la(\a) \neq 0, \a \in L \}$ and $t_j = \la (\sigma_{i_j})$.
 So, $$\la = (\sigma_{i_1}, \cdots , \sigma_{i_1}, \sigma_{i_2}, \cdots , \sigma_{i_2}, \cdots , \sigma_{i_k}, \cdots , \sigma_{i_k})$$
 where $\sigma_{i_j}$ appears $t_j$ times. Note also that $t_1 + t_2
 + \cdots t_k = l-1$. Moreover, for $0\leq a \leq n$, if $a\neq i_j,
 \forall j$, then define $\la^{(a)}$ to be the partition such that
 $\la^{(a)}(\sigma_a) = 1, \la^{(a)}(\sigma_{i_j}) = t_j$; if $a =
 i _v$ for some $v$, then define $\la^{(a)}$ to be the partition such that
 $\la^{(a)}(\sigma_{i_v}) = t_v + 1$, and $\la^{(a)}(\sigma_{i_j}) = t_j, j\neq
 v$. Clearly $\ell(\la^{(a)}) = l$. Observe that for each $\mu \in
 \mathcal{Q}$ with $b_\mu \neq 0$ and $\mu <i> = \la$ for some $i$,
 we have $\mu = \la^{(a)}$ for some $a$.

 Now it is easy to see that the coefficient $c_\la$ of $x_{\la}w'$
 in the representation of the formula $\eqref{h1}$, under the basis
 $\triangle$, is
 \begin{eqnarray*}c_{\la} =
 \sum_{j=0}^nb_{\la^{(j)}}(-s-j)d_j\varphi_{m,j}
 \end{eqnarray*}
 where  $d_{i_v} = t_v +
 1, 1\leq v \leq k$, $d_j = 1$ if $j\neq i_1, i_2, ..., i_k$ and $\varphi_{m,j} = \varphi_{m,j}^{(n,s)}$ that is defined in 2.3.2.
 Meanwhile, one immediately deduces that
 $x_{\la}w'$ does not appear in the representation of the formula
 $\eqref{h2}$. Hence $c_\la = 0$, since $w''$ is a Whittaker vector, that is, for all $m\geq 1$,
\begin{eqnarray}
 \sum_{j=0}^nb_{\la^{(j)}}(-s-j)d_j\varphi_{m,j}=0. \label{co}
 \end{eqnarray}
 Then the assumption that $\varphi$ is a good character implies that
 $ b_{\la^{(j)}} = 0, 0\leq j \leq n.$ But this contradicts with the choice of $\la$.
\QED \vspace{10pt}

\ni 3.2. {\bf Lemma} \ \ {\it Any nontrivial submodule of $V = M/IM$
contains a nonzero Whittaker vector.}

\ni{\bf Proof} \ \  Let $V_1$ be a submodule of $V$. Suppose $V_1$
contains no nonzero Whittaker vector. Use the notation in 2.3.4. Let
$t = max\{ mindeg_1(v) \,|\, v\neq 0, v \in V_1 \}, l = min \{\ell_1
(v) \,|\, mindeg_1 (v) = t, v \in V_1, v\neq 0 \}$. Take a $u \in V_1$
such that $mindeg_1 (u)= t, \ell_1 (u) = l$ (clearly, $l > 0$
). Write $u = \sum\limits_{ 0\geq a \geq t}u_a w'$, where $ u_a =
\sum\limits_{\pi(|\la |) = a}p_{\la, a}x_\la$, with $p_{\la, a}
\in A/I$.  Since $u$ is not a Whittaker vector, there exists a $x
\in \n_\sigma$, for some $\pi(\sigma)
> 0$ such that $u' := xu-\varphi(x)u = \sum\limits_{0 \geq a \geq
t}[x, u_a]w' \neq 0$, where $[x, u_\a]$ stands for
$\sum\limits_{\pi(|\la |) = a}p_{\la, a}[x, x_\la]$. Note that $u'$
is contained in $V_1$. Then, it's easy to see $mindeg_1([x, u_\a]w')
\geq a \geq t$, by Lemma 2.2.3, if $[x, u_\a] \neq 0$. So we have
$mindeg_1(u') \geq t$ and hence $mindeg_1(u') = t $ for the
definition of $t$. In this case, $[x, u_{t}]w' \neq 0$ and
$mindeg_1([x, u_{t}]w') = t$. But this forces $ \ell_1 ([x,
u_{\tau_0}]w' ) < ht(u_t) = \ell_1(u)$ (c.f. Lemma 2.2.3). Thus,
$\ell_1 ( u' ) < l $, which contradicts with the definition of $l$.
\QED \vspace{10pt}

\ni 3.3. {\bf Remark} \ \ 1). Lemma 3.1 and 3.2 suggest that $\B$
satisfies the properties 1) and 2), therefore with the technique
developed in \cite{Wang}, one can set a correspondence of the set of
Whittaker modules and the set of ideals of $S(Z)$. This is treated
in the next section.

2). If $\varphi$ is nonsingular but not a good character, then Lemma
3,1 may not hold. For example, if $\varphi (x_{(2-i, i)}) = 1, i\geq
0$, then $(4x_{\a}^2 + x_{\beta}^2 - 4x_{\a}x_\beta) w'$ is a Whittaker
vector, where $\a = (0, 1), \beta = (-1, 2)$. However, $(4x_\a^2 +
x_\beta^2 - 4x_\a x_\beta) w'$ is not contained in $A/I w'$.

3). Almost all nonsingular characters are good characters. To see
this, denote by $G^{(n,s)}$ the matrix formed by the first $n+1$
rows and columns of $H^{(n.s)}, n\geq 1, s\geq 1$. Then obviously,
the set $\{det(G^{(n,s)})\,|\, n\geq 1, s\geq 1\}$ consists of
countable polynomials of $\varphi (x_{2-a, a}), a\geq 0$. Hence the
statement follows from the fact that a nonsingular character
$\varphi$ is good if all $det(G^{(n,s)})\neq 0, n\geq 1, s\geq 1$.
\vs{18pt}

\cl{\bf\S4. \  Whittaker Modules for $\B$}

The results and their proofs are exactly parallel to \cite{Wang}.
Notation as in 2.1.2. Fix a nonsingular homomorphism $\varphi : \n
\rar \C$, and let $M = M_{\varphi}, w = \bo \otimes 1$. Let $A =
S(Z)$. \vspace{10pt}

4.1.1 {\bf Proposition} \ \ {\it Let $N$ be a submodule of $M =
M_{\varphi}$. Then $N = IM$ for some ideal $I$ of $A = S(Z)$.}

\ni{\bf Proof} \ \ Set $I = \{ x \in A \,|\, xw \in N \}$. One
immediately sees that $I$ is an ideal of $A$ and $IM \subseteq N$.
So we can view $N/IM$ as a submodule of $M/IM$. If $N\neq IM$, then
there exists $pw' \in N/IM$, with $pw' \neq 0, p \in A$, $(w' =
p_I(w))$ by Lemma 3.1 and 3.2. So $pw \in N$ and hence $p \in I$.
Therefore $pw \in IM$, which contradicts with the fact that $pw'
\neq 0$ in $N/IM$. Thus, $N = IM$.\QED

\ni 4.1.2 {\bf proposition} \ \ {\it  Then any nontrivial submodule
of a Whittaker module of type $\varphi$ contains a nontrivial
Whittaker submodule of type $\varphi$.}

\ni{\bf Proof} \ \ It follows immediately from Proposition 4.1.1 and
Lemma 3.2.\QED\vspace{10pt}

\ni4.2. The character $\varphi : \n \rar \C$ naturally extends to a
character of $\uu (\n )$. Let $\uu_{\varphi}(\n )$ be the kernel of
this extension so that $\uu (\n ) = \C \oplus \uu_{\varphi}(\n )$.
Hence,$$ \uu(\B) = \uu(\b_- )\otimes \uu(\n ) = \uu(\b_- ) \oplus
I_{\varphi}$$where $I_{\varphi}=\uu(\B )\uu_{\varphi}(\n )$. For any
$u \in \uu(\B )$, let $u^{\varphi} \in \uu(\b_- )$ be its component
in $\uu(\b_- )$ relative to the above decomposition of $\uu(\B )$.

If $V$ is a Whittaker module generated by a Whittaker vector $v$, let
$\uu_v(\B)$ (resp. $\uu_V(\B )$) be the annihilator of $v$ (resp.
$V$).Then we have, immediately, $V \simeq \uu(\B)/\uu_v(\B)$. Set
$A_V = A \cap \uu_V(\B )$.

 \ni 4.2.1. {\bf Proposition} \ \ {\it Let $V$ be any $\B$ module
that admits a cyclic Whittaker vector $v$. Then $$\uu_v(\B ) =
\uu(\B )A_V + I_\varphi .$$}

{\bf Proof}\ \ Obviously the right hand side of the equation is
contained in the left hand side. So it is enough to show the other
way around. Using the universal property of $M$, we can choose a
surjective homomorphism $\psi : M \rar V$ that sends $w=\bo \otimes
1$ to $v$. Let $Y = ker(\psi )$. Then $Y = IM$ for some $I \subseteq
A$, by 4.1.1.

But for any $x \in \uu_v(\B)$, i.e. $xv=0$, one has $x^\varphi v=0$ and
hence $x^\varphi w \in Y$. Then $x^\varphi w = \sum_ip_ix_iw,\, x_i
\in \uu'', p_i \in I$. Thus, $x^\varphi = \sum\limits_ip_ix_i
\subset \uu I$. But clearly $I \subseteq A_V$, therefore $x^\varphi
\in \uu A_V$. So, $x \in \uu A_V + I_\varphi$. \QED \vspace{6pt}

\ni4.2.2. {\bf Theorem} \ \ {\it  The correspondence $$V \rar
A_V$$sets up a bijection between the set of all the isomorphic
classes of Whittaker modules for $\B$ and the set of all the ideals
of $A = S(Z)$.}

\ni{\bf Proof} \ \ Note that for any $I \subseteq A$, if let $V =
M/IM$, then $A_V = I$. Now the theorem follows
immediately.\vspace{6pt}

\ni 4.2.3. {\bf Corollary} \ \ {\it  For any maximal ideal
$\mathfrak{m} \in S(Z)$, $L_{\varphi, \mathfrak{m}} =
M/\mathfrak{m}M$ is simple and any simple Whittaker module of type
$\varphi$ is of this form.}

\ni{\bf Proof} \ \ Observe that if $I \subseteq J \subseteq S(Z)$,
then $M/JM$ is a quotient of $M/IM$. The corollary now follows from
Theorem 4.2.2 immediately.\QED

\end{document}